\documentclass[12pt,a4paper]{article}

\usepackage{amsfonts,amsmath,amssymb,amsopn}

\title{The filtered Poincar\'e lemma in higher level \\ (with applications to algebraic groups)}

\author{Bernard Le Stum and Adolfo Quir\'os
 \thanks{Our collaboration was made possible by the MCRTN Arithmetic Algebraic Geometry (MRTN-CT-2003-504917) of the European Union.
The second author was partially supported by MCYT (Spain) project BFM2003-02606.}}

\date{}

\newcommand{\Char}{\mathop{\rm char}\nolimits}
\newcommand{\Id}{\mathop{\rm Id}\nolimits}
\newcommand{\Hom}{\mathop{\rm Hom}\nolimits}

\newcommand{\Spec}{\mathop{\rm Spec}\nolimits}
\newcommand{\Ob}{\mathop{\rm Ob}\nolimits}
\newcommand{\Arr}{\mathop{\rm Arr}\nolimits}
\newcommand{\Fil}{\mathop{\rm Fil}\nolimits}
\newcommand{\Sch}{\mathop{\rm Sch}\nolimits}
\newcommand{\Cris}{\mathop{\rm Cris}\nolimits}
\newcommand{\cris}{\mathop{\rm cris}\nolimits}
\newcommand{\CRIS}{\mathop{\rm CRIS}\nolimits}
\newcommand{\ZAR}{\mathop{\rm ZAR}\nolimits}

\def \Z {\mathbf Z}
\def \C {\mathbf C}

\def \N {\mathbf N}

\begin{document}

\maketitle

\begin{abstract}

We show that the Poincar\'e lemma we proved elsewhere in the context of crystalline cohomology of higher level behaves
well with regard to the Hodge filtration. This allows us to prove the Poincar\'e lemma for transversal crystals of
level $m$. We interpret the de Rham complex in terms of what we call the Berthelot-Lieberman construction, and show how
the same construction can be used to study the conormal complex and invariant differential forms of higher level for a
group scheme. Bringing together both instances of the construction, we show that crystalline extensions of transversal
crystals by algebraic groups can be computed by reduction to the filtered de Rham complexes. Our theory does not ignore
torsion and, unlike in the classical ($m=0$), not all closed forms are invariant. Therefore, close invariant
differential forms of level $m$ provide new invariants and we exhibit some examples as applications.

\end{abstract}

\tableofcontents

%
%
\section*{Introduction}

In a series of articles, starting with \cite{LeStumQuiros97} and \cite{LeStumQuiros01}, we are using the partial divided powers of Berthelot to study the geometry of algebraic varieties of positive characteristic.
This gives new insight into the p-adic cohomological theories.
Unlike other works on the subject (\cite{EtesseLeStum97}, \cite{Trihan98} and \cite{Crew03*}), we do not use crystalline cohomology of higher level as a tool to obtain results in rigid cohomology and, in particular, we do not ignore torsion.
In fact, torsion is very rich in this theory and provides new invariants that help understand the geometry of algebraic varieties.
For example, we will show how the sheaf of closed invariant differential forms of higher level can tell you exactly where the supersingular locus of a family of elliptic curves is.

\bigskip

Following A. Ogus in \cite{Ogus94}, we introduced in \cite{LeStumQuiros97} the notion of transversal crystal of higher level.
Although we could slightly improve on some of Ogus' results, progress was hampered by the lack of de Rham techniques for computing exactly crystalline cohomology in higher level.
A natural answer to this problem was provided in \cite{LeStumQuiros01}, where we developed ideas of P. Berthelot, introducing the de Rham complex of higher level and proving the exact Poincar\'e lemma.
In the present article, we extend the Poincar\'e lemma of higher level to transversal crystals by paying close attention to filtrations.
This is used to give a precise description of the group of extensions of a transversal crystal by a commutative group scheme.
In a forthcoming article, we want to use these results to show that Dieudonn\'e crystals of higher level are transversal.

\bigskip

After reviewing in section 1 a few results on filtrations, and especially fixing terminology about filtered derived categories, in section 2 we verify that the formal Poincar\'e lemma in higher level behaves as expected with respect to filtrations.
In doing so, we reinterpret the de Rham and linearized de Rham complexes introduced in \cite{LeStumQuiros01} as particular cases of what we call a Berthelot-Lieberman complex.
This general construction will be used in section 4 to define the conormal complex of higher level of a group scheme.
But, before, we prove the filtered Poincar\'e lemma for transversal crystals.
This is done in section 3 and requires careful attention to the behavior of the filtrations all along the process.
In section 4 which is completely independent of the preceding one, we study the relation between the conormal sheaf of higher level and invariant differential forms of higher level.
Unlike in the classical case (level 0), not all invariant forms are closed.
Actually, the module of closed invariant forms is isomorphic to the first cohomology group of the conormal complex.
We present concrete examples, including the Legendre family of elliptic curves, and give the relation with de Rham cohomology of higher level in the case of abelian schemes.
Section 5 brings together the two previous sections.
More precisely, we show that crystalline extensions of transversal crystals by algebraic groups can be computed by reduction to filtered de Rham complexes.
As an application, we show that the extension group of the partial divided power ideal by a smooth group is nothing else but a lifting of the module of closed invariant differentials of higher level.

\bigskip

Most results here are inspired by theorems that have been well known for a long time in the case of usual divided
powers and classical crystalline cohomology. In particular, this work owes much to P. Berthelot, L. Breen, L. Illusie,
W. Messing and A. Ogus.

%
%
\section*{Conventions}

Starting at section \ref{secformal}, we let $p$ be a prime, $m \in \N$ and, unless $m = 0$, all schemes are assumed to be $\Z_{(p)}$-schemes.

%
%
\section{Generalities about filtrations}

Concerning filtrations, we use the terminology of \cite {Deligne71b}, 1.1.
In particular, we will only consider filtrations of type $\Z$ in abelian categories.
Also, filtration will always mean decreasing filtration.

\subsection{Definition} A filtration $\Fil^{\bullet}$ on an object $M$ is \emph{effective} if $\Fil^0M = M$.

\bigskip

Unless otherwise specified, we will only consider effective decreasing filtrations.
Note that the filtration induced on a subquotient (cf. \cite {Deligne71b}, 1.1.10) by an effective filtration is still effective.
One can also check that the image of an effective filtration by a semi-exact (meaning left or right exact) multiadditive functor (cf. \cite {Deligne71b}, 1.1.12) is still effective.
This applies in particular to the tensor product filtration.

\subsection{Definition} The \emph{trivial (effective) filtration} on an object $M$ is given by
$$
\Fil^kM = \left \{ \begin{array} {cc} M \mbox{ if } k \leq 0 \\ 0 \mbox{ if } k > 0\end{array} \right..
$$

\bigskip

In \cite {Deligne71b} (definition 1.3.6), filtered quasi-isomorphisms are only defined for so-called biregular filtrations.
This definition does not generalize well.
As in \cite {Ogus94}, section 4.4, we will need a more restrictive notion.

\subsection{Definition} A morphism of filtered complexes $M^\bullet \to N^\bullet$ is a \emph{true filtered quasi-isomorphism} if, for each $k \in \Z$, the induced morphism $\Fil^kM^\bullet \to \Fil^kN^\bullet$ is a quasi-isomorphism.
A \emph{filtered homotopy} is a morphism of filtered complexes which is a homotopy.
Two morphisms of filtered complexes $f, g : M^\bullet \to N^\bullet$ are \emph{homotopic} if there exists a filtered homotopy $h : M^\bullet \to N^\bullet$ such that $g$ is homotopic to $f$ with respect to $h$.

\bigskip

As explained in \cite {Ogus94}, page 80, the above notions of filtered homotopy and true  filtered quasi-isomorphisms are suitable to define the filtered derived category of a Grothendieck category, e.g. a category of modules.
Moreover, any left exact additive functor $F$ between such categories gives rise to a filtered derived functor and this construction is completely compatible with the non-filtered situation in the sense that we always have
$$
RF(\Fil^kM^\bullet) = \Fil^kRFM^\bullet.
$$
In particular, there exists a canonical spectral sequence
$$
E_1^{i,j}  = R^{i+j}F(Gr^iM^\bullet) \Rightarrow R^nF(M^\bullet)
$$
that endows, for each $n \in \Z$, $R^nF(M^\bullet)$ with a canonical filtration.

\bigskip

We now need to recall some definitions and results from \cite {Ogus94} and \cite {LeStumQuiros97} on transversal filtrations.

\subsection{Definition} A \emph{ring filtration} on a ring $A$ is an effective filtration by ideals $I^{(k)}$ such that $I^{(k)}I^{(l)} \subset I^{(k+l)}$.
A filtered $A$-module $(M, \Fil^\bullet)$ is \emph{transversal} (resp. \emph{almost transversal}) to $I^{(\bullet)}$ if for each $k \in \Z$,
$$
I^{(1)}M \cap \Fil^kM = (resp. \subset) \sum_{i + j = k, i > 0} I^{(i)}\Fil^{j}.
$$
If $(M, \Fil^\bullet)$ is a filtered $A$-module, the \emph{saturation} of the filtration with respect to $I^{(\bullet)}$ is the tensor product filtration  $\overline {\Fil}^\bullet$ on $M$ under the identification $A \otimes_A M = M$.

\bigskip

When $A$ is filtered by the powers of an ideal $I$, we simply say transversal (resp. almost transversal, resp.  saturation with respect) to $I$.
Note that a module filtration is always transversal to $0$ and that, if a filtration is almost transversal, then its saturation is transversal.

\subsection{Definition} The \emph{trivial transversal filtration} on an $A$-module $M$ is the saturation of the trivial filtration.

\bigskip

Note that the trivial transversal filtration is given by
$$
\overline{\Fil}^kM = \left \{ \begin{array} {cc} M \mbox{ if } k \leq 0 \\ I^{(k)}M \mbox{ if } k > 0\end{array} \right.
$$
%
%
\section{The formal filtered Poincar\'e lemma} \label {secformal}

The aim of this section is to give a filtered version of theorem 3.3 of \cite{LeStumQuiros01}.
We first recall some properties of $m$-PD-envelopes and introduce the notion of Berthelot-Lieberman complex.

\bigskip

If $X \hookrightarrow Y$ is an immersion of schemes, then we will denote by $P_{Xm}(Y)$ its divided power envelope of level $m$.
We will write $\mathcal P_{Xm}(Y)$ for the structural sheaf and $\mathcal I_{Xm}(Y)$ the $m$-PD-ideal.
 We will use a superscript $n$ to denote the same objects modulo $\mathcal I^{\{n + 1\}}_{Xm}(Y)$.
If $X$ is an $S$-scheme, we will denote by $P_{X/Sm}(r)$ the $m$-PD envelope of the diagonal embedding of $X$ in $X^{r+1}$ and modify all other notations accordingly.

\bigskip

The notion of $m$-PD envelope is functorial in the sense that any commutative diagram
$$
\begin{array}{lcl}
X' & \hookrightarrow &  Y' \\
 \downarrow g     &   &     \downarrow g'\\
X &  \hookrightarrow & Y
\end{array}
$$
canonicaly extends to
$$
\begin{array}{ccccc}
X' &  \hookrightarrow & P_{X'm}(Y') &   \hookrightarrow & Y' \\
\downarrow     & & \downarrow     & &  \downarrow \\
X & \hookrightarrow & P_{Xm}(Y) &   \hookrightarrow & Y
\end{array}
$$

We recall the following fundamental results of Berthelot :

\subsection{Proposition} \label{gflat} \emph{
Let $ f : X \hookrightarrow Y$ be an immersion of schemes. Then,
\begin{enumerate}
\item  If $f$ has locally a smooth retraction then, for all $n$, $\mathcal P^n_{Xm}(Y)$ is a locally free $\mathcal O_X$-module of finite rank. Actuallly, $\mathcal P_{Xm}(Y)$ itself is locally free when $Y$ is uniformly killed by a power of $p$.
\item Assume that $f$ has a smooth retraction $p : Y \to X$ and let
$$
\begin{array}{lcl}
X' & \leftarrow  \atop\hookrightarrow &  Y' \\
 \downarrow g     &   &     \downarrow g'\\
X & \leftarrow \atop \hookrightarrow & Y
\end{array}
$$
be a be cartesian diagram  (both ways).
Then, we have for all $n$,
$$
g^*\mathcal P^n_{Xm}(Y) \simeq \mathcal P^n_{X'm}(Y').
$$
and even
$$
g^*\mathcal P_{Xm}(Y) \simeq \mathcal P_{X'm}(Y').
$$
if $Y$ is uniformly killed by a power of $p$.
\end{enumerate}}

\bigskip

{\bf Proof : } This follows from propositions 1.4.6 and 1.5.3 of \cite {Berthelot96}.

\bigskip

Let $X \hookrightarrow Y(\bullet)$ be an immersion of a constant scheme into a simplicial complex.
In other words, we are given a family of immersions of $X \hookrightarrow Y(r)$ compatible with the differentials $d_i$ and the degeneracy arrows $s_i$ of $Y(\bullet)$.
Taking $m$-PD envelopes gives rise to a simplicial complex $P_{Xm}(Y\bullet)$ from which we derive a complex $(\mathcal P_{Xm}(Y\bullet), \mathcal I^{\{k\}}_{Xm}(Y\bullet))$ of filtered rings.
We then consider as in \cite{Illusie71} the normalization of $\mathcal P_{Xm}(Y\bullet)$ which is the subcomplex of ideals defined by $N \mathcal P_{Xm}(Yr) := \cap \ker s_i^*$.
Finally, we are interested in the quotient of $N \mathcal P_{Xm}(Y\bullet)$ by the differential subalgebra generated by the ideal $\mathcal I^{\{p^m+1\}}_{Xm}(Y1)$ of  $\mathcal P_{Xm}(Y1)$, which we write $\Omega^\bullet_{Xm}(Y)$.

\subsection{Definition} The complex $\Omega^\bullet_{Xm}(Y)$ is the \emph{Berthelot-Lieberman complex} of $X \hookrightarrow Y(\bullet)$.
The \emph{Hodge filtration} on $\Omega^\bullet_{Xm}(Y)$ is the filtration $\Fil^\bullet_H$ induced by the filtration of $\mathcal P_{Xm}(Y\bullet)$.

\subsection{Definition} \begin{enumerate}
\item In a category with products, the \emph{product simplicial complex} $X^{prod}(\bullet)$ of an object $X$ is defined by $X^{prod}(r) = X^{r+1}$,
$$
d_i(x_1, \ldots, x_{r+2}) = (x_1, \ldots, x_i, x_{i+2}, \ldots, x_{r+2})
$$
and
$$
s_i(x_1, \ldots, x_{r}) = (x_1, \ldots, x_i, x_{i+1}, x_{i+1}, x_{i+2}, \ldots, x_{r+2}).
$$
\item If $X(\bullet)$ is a simplicial complex in any category, the \emph{shifted simplicial complex} is defined by $X^+(r) = X(r+1), d^+_i = d_{i+1}, s^+_i = s_{i+1}$.
\end{enumerate}

\bigskip

We can now reformulate some definitions from section 1 of \cite {LeStumQuiros01}.

\bigskip

Let $S$ be a scheme with $p$ locally nilpotent and $X$ an $S$-scheme.

\subsection{Definition} The Berthelot-Lieberman complex $\Omega^\bullet_{X/Sm}$ of the diagonal embedding $X \hookrightarrow (X/S)^{prod}(\bullet)$, is the \emph{de Rham complex of level $m$}.
The Berthelot-Lieberman complex $L_X(\Omega^\bullet_{X/Sm})$ of the shifted simplicial complex $(X/S)^{prod+}(\bullet)$ is the \emph{linearized de Rham complex  of level $m$}.

\bigskip

Note that when $X$ is smooth, the de Rham complex of level $0$ is the usual de Rham complex with its Hodge filtration.
Note also that $L_X(\Omega^\bullet_{X/Sm})$ is what we called $\Omega^\bullet_{P/Sm}$ in \cite {LeStumQuiros01}.

\bigskip

For the rest of this section, we assume that $X$ is a smooth scheme over $S$.

\subsection{Remarks} \label{local} Recall from section 1 of \cite {LeStumQuiros01} that, if we have local coordinates $t_1, \ldots, t_n$ on $X$ and, as usual, we set $\tau_i := 1 \otimes t_i - t_i \otimes 1$, then $\mathcal P_{X/Sm}(r)$ is a free $\mathcal O_X$-module on generators $\tau^{\{J_1\}} \otimes \cdots \otimes \tau^{\{J_r\}}$ with $|J_i| \geq 0$ (we use the standard multiindex convention).
Using proposition 1.5.3 of \cite {Berthelot96}, one easily checks that $\mathcal I_{X/Sm}^{\{k\}}(r)$ is generated by the $\tau^{\{J_1\}} \otimes \cdots \otimes \tau^{\{J_1\}}$ with $|J_1| + \cdots |J_r| \geq k$.

\bigskip

Now, if as usual again, $dt_i$ denotes the image of $\tau_i$ in $\Omega^1_{X/Sm}$, we know that $\Omega^1_{X/Sm}$ is a free $\mathcal O_X$-module on generators $(dt)^J$ with $0 < |J| \leq p^m$.
As we discussed in section 1 of \cite{LeStumQuiros01}, even if $\Omega^r_{X/Sm}$ is generated by the $(dt)^{J_1} \otimes \cdots \otimes (dt)^{J_r}$ with $0 < |J_i| \leq p^m$, there are some relations among these generators.
Anyway, we see that $\Fil^k_H \Omega^r_{X/Sm}$ is generated by the $(dt)^{J_1} \otimes \cdots \otimes (dt)^{J_1}$ with $0 < |J_i| \leq p^m$ and $|J_1| + \cdots |J_r|  \geq k$.
In particular, $\Fil^k_H = 0$ for $k > rp^m$.

\bigskip

Finally $L_X(\Omega^r_{X/Sm})$ is generated as $\mathcal P_{X/Sm}(1)$-module by the $(d\tau)^{J_1} \otimes \cdots \otimes (d\tau)^{J_r}$ \linebreak with $0 < |J_i| \leq p^m$, and $\Fil^k_H L_X(\Omega^r_{X/Sm})$ is the $\mathcal O_X$-submodule generated by the \linebreak $\tau^{\{I\}} (d\tau)^{J_1} \otimes \cdots \otimes (d\tau)^{J_r}$ with $|I| + |J_1| + \cdots |J_r| \geq k$.

\bigskip

A quick look at the proof of proposition 1.4 in \cite{LeStumQuiros01} shows that the Hodge filtration is a filtration by locally free submodules.
In fact, this question is local and we may then proceed by extracting a basis from a given set of generators.
But all the relations are homogeneous and therefore, respect the Hodge filtration.

\bigskip

Recall that if $\mathcal F$ is any $\mathcal D_{X/S}^{(m)}$-module, we can form the de Rham complex
$$
\mathcal F \otimes_{\mathcal O_X} \Omega^\bullet_{X/Sm}.
$$
of $\mathcal F$.
In particular, we have the \emph{$m$-connection} $\mathcal F \to \mathcal F \otimes_{\mathcal O_X} \Omega^1_{X/Sm}.$
For further use, recall also that the de relative de Rham cohomology of level $m$ of $\mathcal F$ is
$$
\mathcal H^n_{dRm}(X, \mathcal F) := \mathbf R^np_{X*}\mathcal F \otimes_{\mathcal O_X} \Omega^\bullet_{X/Sm}
$$
where $p_{X} : X \to S$ si the structural map.

\bigskip

Concerning Griffiths tranversality, we refer to section 2.2 of \cite {LeStumQuiros97}.
Let us just recall that a filtration $\Fil^\bullet$ on a $\mathcal D_{X/S}^{(m)}$-module $\mathcal F$ is a filtration by $\mathcal O_X$-submodules and that it is said \emph{Griffiths tranvsersal} if we always have
$$
\mathcal D_{X/Sj}^{(m)} \Fil^k \mathcal F \subset \Fil^{k - j} \mathcal F.
$$

\subsection{Proposition}\^{E}\emph{A filtered $\mathcal D_{X/S}^{(m)}$-module is Griffiths transversal if and only if the  $m$-connection is compatible with the filtrations.
\bigskip}

{\bf Proof : } This is a local question and, if we have local coordinates, the $m$-connection on $\mathcal F$ is given by
$$
s \mapsto \sum_{0 < |J| \leq p^m} \partial^{[J]}s \otimes (dt)^J.
$$
If the filtration is Griffiths transversal and $s \in \Fil^k\mathcal F$, we have
$$
\partial^{[J]}s \otimes (dt)^J \in \Fil^{k-|J|}\mathcal F \otimes \Fil^{|J|}\Omega^\bullet_{X/Sm} \subset \Fil^k [\mathcal F \otimes \Omega^\bullet_{X/Sm}]
$$
It follows that the connection preserves the filtration.
Conversely, if the $m$-connection is compatible with the filtration and $s \in \Fil^k\mathcal F$, then, necessarily, for all $i = 1, \ldots, n$ and $j \leq m$, we have
$$
\partial^{[p^j]}s \otimes (dt)^{p^j} \in  \Fil^k [\mathcal F \otimes \Omega^\bullet_{X/Sm}]
$$
so that, necessarily, $\partial_i^{[p^j]}s \in \Fil^{k-|p^j|}\mathcal F$ and by remark 2.2.2 (ii) of \cite{LeStumQuiros97}, this is sufficient for Griffiths tranversality.

\bigskip

Since the differential on a de Rham complex can be described by the Leibnitz rule, we have the following :

\subsection{Corollary} \emph{If $\mathcal F$ is a filtered $\mathcal D_{X/S}^{(m)}$-module and if, for each $r$, $\mathcal F \otimes \Omega^r_{X/Sm}$ is endowed with the tensor product filtration, then $\mathcal F$ is Griffiths transversal if and only if the de Rham complex of $\mathcal F$ is a filtered complex.}

\subsection{Lemma} \label{inverse} \emph{For all $r$, if we endow $\Omega^r_{X/Sm}$ and $L_X(\Omega^r_{X/Sm})$ with their Hodge filtration and $\mathcal P_{X/Sm}(1)$ with its $m$-PD-filtration, we have an isomorphism of filtered modules
$$
L_X(\Omega^r_{X/Sm}) = \mathcal P_{X/Sm}(1) \otimes_{\mathcal O_X} \Omega^r_{X/Sm}.
$$
Moreover, the filtration on $\mathcal P_{X/Sm}(1) \otimes_{\mathcal O_X} \Omega^r_{X/Sm}$ is the saturation with respect to $\mathcal I_{X/Sm}^{\{\bullet\}}(1)$ of the inverse image by $\mathcal P_{X/Sm}(1) \to \mathcal O_X$ of the Hodge filtration on $\Omega^r_{X/Sm}$.}

\bigskip

{\bf Proof : } Since we assume that $X$ is smooth over $S$, then as mentioned in section 1.5 of \cite{LeStumQuiros01}, $\mathcal P^+_{X/Sm}(r)$ is canonically isomorphic as filtered algebra to $\mathcal P_{X/Sm}(1) \otimes_{\mathcal O_X} \mathcal P_{X/Sm}(r)$.
The assertion then follows from the functoriality of our construction.
The second assertion is local in nature and follows directly from the above local description of our filtrations.

\bigskip

\subsection{Proposition} \emph{The Hodge filtration on $L_X(\Omega^\bullet_{X/Sm})$ is transversal to $\mathcal I_{X/Sm}^{\{\bullet\}}(1)$.}

\bigskip

{\bf Proof : } The morphism of filtered rings $(\mathcal P_{X/Sm}(1),  \mathcal I_{X/Sm}^{\{\bullet\}}(1)) \to (\mathcal O_X, 0)$ obviously satisfies the assumptions of proposition 1.1.8 of \cite {LeStumQuiros97}.
Since any filtration is transversal to the $0$-ideal, we see that the inverse image of the Hodge filtration is almost transversal and it follows that its saturation is transversal.
Our assertion now results from lemma \ref{inverse}

\bigskip

The following is a generalization of theorem 3.3 of \cite {LeStumQuiros01}

\subsection {Proposition} \label{formal}  \emph{If $\mathcal O_X$ is endowed with the trivial filtration and  $L_X(\Omega^\bullet_{X/Sm})$ with the Hodge filtration, the canonical map
$$
\mathcal O_X \to L_X(\Omega^\bullet_{X/Sm})
$$
is a true filtered quasi-isomorphism.
More precisely, locally on $X$, it is a filtered homotopy equivalence.}

\bigskip

{\bf Proof : } The proof works exactly as in \cite {LeStumQuiros01} once one notices that the homotopy of 2.1 in \cite {LeStumQuiros01} is a filtered homotopy.

%
%
\section{The filtered Poincar\'e lemma}

We will explain here how the results of \cite{LeStumQuiros01}, section 4 extend to the case of transversal $m$-crystals.

\bigskip

Let $(S, \mathfrak a, \mathfrak b)$ be a $m$-PD-scheme with $p$ locally nilpotent and $p \in \mathfrak a$.
Let $X$ be an $S$-scheme to which the $m$-PD-structure of $S$ extends.

\bigskip

Unless otherwise specified, we will assume in this section that $X$ is smooth over $S$.

\subsection{Notations} If $Y$ is any object in a topos $\mathcal T$, we will denote by $\mathcal T_{/Y}$ the localized category and by $j_Y : \mathcal T_{/Y} \to \mathcal T$ the restriction map.

\bigskip

In particular, we will consider here the restriction map
$$
j_X : (X/S)^{(m)}_{\cris|X} \to (X/S)^{(m)}_{\cris}
$$
and the projection
$$
u_X : (X/S)^{(m)}_{\cris} \to X_{Zar}.
$$

\subsection{Proposition} \label {ideals} \emph{The ideal $\mathcal K^{(m)}_{X/S} := j_{X*}j_X^{-1} \mathcal I^{(m)}_{X/S}$ of $ j_{X*}j_X^{-1} \mathcal O^{(m)}_{X/S}$ is an $m$-PD-ideal, there is a canonical exact sequence
$$
0 \longrightarrow  \mathcal K^{(m)}_{X/S} \longrightarrow  j_{X*}j_X^{-1} \mathcal O^{(m)}_{X/S} \longrightarrow u_X^{-1} \mathcal O_X \longrightarrow 0
$$
and for all $k \in \Z$, we have ${\mathcal K^{(m)}_{X/S}}^{\{k\}} = j_{X*}j_X^{-1} [{\mathcal I^{(m)}_{X/S}}^{\{k\}}]$.}

\bigskip

{\bf Proof : } One easily checks that, if $\mathcal F$ is any $\mathcal O_X$-module, the adjunction map
$$
j_{X*}j^{-1}_Xu_X^{-1} \mathcal F \to u_X^{-1} \mathcal F
$$
is an isomorphism.
In particular, if we apply the exact functor $j_{X*}j^{-1}_X$ to the exact sequence
$$
0 \longrightarrow  \mathcal I^{(m)}_{X/S} \longrightarrow \mathcal O^{(m)}_{X/S} \longrightarrow u_X^{-1} \mathcal O_X \longrightarrow 0,
$$
we get the expected exact sequence
$$
0 \longrightarrow  j_{X*}j^{-1}_X \mathcal I^{(m)}_{X/S} \longrightarrow  j_{X*}j_X^{-1} \mathcal O^{(m)}_{X/S} \longrightarrow u_X^{-1} \mathcal O_X \longrightarrow 0.
$$
Thus, we see that if $U \subset X$ is any open subset and $U \hookrightarrow Y$ an $m$-PD-thickening, we have the following exact sequence
$$
0 \longrightarrow (K^{(m)}_{X/S})_Y \longrightarrow \mathcal O_P \longrightarrow \mathcal O_U \longrightarrow 0
$$
where $P := P_{Um}(Y \times_S X)$.
Hence, we see that $(K^{(m)}_{X/S})_Y$ is the $m$-PD-ideal $\mathcal I_P$ of $P$.
It follows that $K^{(m)}_{X/S}$ is an $m$-PD-ideal and that for all $k \in \Z$, $({K^{(m)}_{X/S}}^{\{k\}})_Y = \mathcal I^{\{k\}}_P = ({\mathcal I^{(m)}_{X/S}}^{\{k\}})_Y.$

\bigskip

Recall that $u_{|X} := u_X \circ j_X$ is in a natural way a morphism of ringed sites.

\subsection{Definition} If $\mathcal F$ is a filtered $\mathcal O_X$-module, its \emph{linearization (of level $m$)} is $L^{(m)}(\mathcal F) := j_{X*}u_{|X}^*\mathcal F$  endowed with the saturation of the filtration $L^{(m)}(\mathcal \Fil^k\mathcal F)$ with respect to the $m$-PD-ideal ${\mathcal K^{(m)}_{X/S}}^{\{\bullet\}}$ of $L(\mathcal O^{(m)}_{X/S})$.

\subsection{Lemma} \emph{If $\mathcal F$ is a filtered $\mathcal O_X$-module,we have a canonical isomorphism of filtered modules
$$
L^{(m)}(\mathcal F)_X = \mathcal P_{X/Sm}(1) \otimes_{\mathcal O_X} \mathcal F.
$$}

{\bf Proof : } Thanks to proposition 4.3 (1) of \cite{LeStumQuiros01}, only the assertion concerning the filtration has to be checked.
Since saturation is just tensor product with the ideal filtration, it is sufficient to note that, as in the proof of lemma \ref{ideals}, $(\mathcal K^{(m)}_{X/S})_X := \mathcal I_{X/Sm}(1)$.

\bigskip

We recall that $L^{(m)}(\Omega^\bullet_{X/Sm})$ has a natural structure of complex of crystals whose realization on $X$ is nothing but the linearized de Rham complex $L_X(\Omega^\bullet_{X/Sm})$.
We can now state the formal Poincar\'e lemma in its crystalline form.

\subsection{Theorem} \label {crysform} \emph{If $E$ is a filtered $\mathcal O^{(m)}_{X/S}$-module which is saturated with respect to $\mathcal I^{(m)}_{X/S}$, then the morphism
$$
E \to E \otimes_{\mathcal O^{(m)}_{X/S}} L^{(m)}(\Omega^\bullet_{X/Sm})
$$
is a true filtered quasi-isomorphism.
More precisely, locally on $\Cris^{(m)}(X/S)$, it is a filtered homotopy equivalence.}

\bigskip

{\bf Proof : } Since the filtration on $E$ is saturated, the canonical identification
$$
E \otimes_{\mathcal O^{(m)}_{X/S}} \mathcal O^{(m)}_{X/S} = E
$$
is compatible with the filtrations.
Since the pull-back of a filtered homotopy is still filtered homotopy, the standard arguments allow us to reduce to the case $E = \mathcal O^{(m)}_{X/S}$ and then to theorem \ref{formal}

\bigskip

The definition of a transversal $m$-crystal is given in \cite {LeStumQuiros97}, section 4.
Roughly speaking, it is a crystal of transversal modules, but the reader should consider looking at the above reference if he really wants a precise definition as well as a description of its relation with Griffiths tranversality.

\subsection{Proposition} \emph{Let $E$ be a transversal $m$-crystal on $X/S$.
If $\mathcal F$ is any filtered $\mathcal O_X$-module, there is a canonical isomorphism of filtered modules
$$
E \otimes L^{(m)}(\mathcal F) \simeq L^{(m)}(E_X \otimes \mathcal F).
$$}

{\bf Proof : } We already have an isomorphism of crystals thanks to proposition 4.3 (3) of \cite {LeStumQuiros01}.
More precisely, if If $U \subset X$ is any open subset and $U \hookrightarrow Y$ an $m$-PD-thickening, we have an isomorphism
$$
E_Y \otimes \mathcal O_P \otimes \mathcal F \simeq E_P \otimes \mathcal F \simeq \mathcal O_P \otimes E_Y \otimes \mathcal F,
$$
where $P := P_{Um}(Y \times_S X)$ as in the proof of proposition \ref {ideals}, and we need to show that it is compatible with filtrations.
But this follows from the definition in \cite{LeStumQuiros97}, 4.2.2 of a transversal $m$-crystal.

\subsection{Corollary} \label{cryslem} \emph{If $E$ is a transversal $m$-crystal on $X/S$, there is a canonical true filtered quasi-isomorphism
$$
E \to L^{(m)}(E_X \otimes_{\mathcal O_X} \Omega^\bullet_{X/Sm}).
$$
More precisely, locally on $\Cris^{(m)}(X/S)$, it is a filtered homotopy equivalence.}

\bigskip

In order to obtain the filtered Poincar\'e lemma at the cohomological level, we need the following :

\subsection{Proposition} \label {rustar} \emph{If $\mathcal F$ is a filtered $\mathcal O_X$-module, we have a canonical isomorphism in the filtered derived category $Ru_{X*}L^{(m)}(\mathcal F) = \mathcal F$.}

\bigskip

{\bf Proof : } We must show that, for all $k \in \Z$, we have $R^iu_{X*}\Fil^kL^{(m)}(\mathcal F) = 0$ for $i > 0$ and $u_{X*}\Fil^kL^{(m)}(\mathcal F) = \Fil^k\mathcal F$.
Using proposition \ref{ideals},
$$
\Fil^kL^{(m)}(\mathcal F) = \sum_{i + j = k, i \geq 0} (j_{X*} j_X^{-1} {\mathcal I^{(m)}_{X/S}}^{\{i\}} )(j_{X*}u_{|X}^*(\Fil^{j}\mathcal F)).
$$
Since $j_{X*}$ is exact, it follows that $\Fil^kL^{(m)}(\mathcal F) = j_{X*} E$ with
$$
E :=  \sum_{i + j = k, i \geq 0} j_X^{-1} {\mathcal I^{(m)}_{X/S}}^{\{i\}}u_{|X}^*(\Fil^{j}\mathcal F).
$$
Regarding the higher direct images, it is sufficient to recall, that since $j_{X*}$ and $u_{|X*}$ are exact, it is automatic that $R^iu_{X*} j_{X*} E = 0$ for $i > 0$.

Now, we know from  proposition 4.3 (2) of \cite {LeStumQuiros01} that, for any $\mathcal O_X$-module $\mathcal F$, we have \linebreak $u_{X*}L^{(m)}(\mathcal F) = \mathcal F$.
Also, $u_X^{-1}$ being fully faithful, $u_{X*}u_X^{-1} \mathcal F = \mathcal F$.
Since $u_{X*}$ is exact and $\mathcal K^{(m)}_{X/S}$ is the kernel of the natural map $L^{(m)}(\mathcal O_X) \to u_X^{-1} \mathcal O_X$, we see that $u_{X*}\mathcal K^{(m)}_{X/S} = 0$.
It follows that $u_{X*}$ ignores saturation with respect to $\mathcal K^{(m)}_{X/S} $ and so,
$$
u_{X*}\Fil^kL^{(m)}(\mathcal F) = u_{X*}L^{(m)}(\Fil^k\mathcal F) = \Fil^k\mathcal F
$$
as asserted.

\subsection{Theorem} \emph{If $E$ is a transversal $m$-crystal on $X/S$, there is a canonical isomorphism in the filtered derived category
$$
Ru_{X*} E \simeq E_X \otimes \Omega^\bullet_{X/Sm}.
$$}

\bigskip

{\bf Proof : } Using proposition \ref{rustar}, this is an direct consequence of \ref{cryslem}.

\subsection {Corollary} \emph{Even if we no longer assume $X$ smooth, but if $i : X \hookrightarrow Y$ is an embedding into a smooth $S$-scheme and if $E$ is a transversal $m$-crystal on $X$, then there is a canonical isomorphism in the filtered derived category
$$
i_*Ru_{X*} E \simeq (i_{\cris*}E)_Y \otimes \Omega^\bullet_{Y/Sm}.
$$}

\bigskip

{\bf Proof : } Works exactly as at the end of section 4 in \cite {LeStumQuiros01}.

\subsection {Corollary} \emph{In the situation of the previous corollary, there is a canonical filtered isomorphism
$$
R\Gamma((X/S)^{(m)}_{\cris}, E) \simeq R\Gamma(Y, (i_{\cris,m*}E)_Y \otimes \Omega^\bullet_{Y/Sm})).
$$
In other words, for all $i$, we have
$$
H^i_{\cris,m}(X, \Fil^k E) = H^i(\Fil^k((i_{\cris,m*}E)_Y \otimes \Omega^\bullet_{Y/Sm})).
$$}

%
%

\section{Differentials of higher level on a group scheme}

In this section, which is independent of the previous one, we define the conormal complex of level $m$ of a group scheme and study invariant differential forms of higher level.

\subsection{Definition} The \emph{simplicial complex $G^{gr}(\bullet)$ associated to a group $G$} in a category with products has components $G^{gr}(r) := G^{r}$, differentials $d_i : G^{r+1} \to G^{r}$ defined by
$$
\begin{array}{lcl }
 d_0(g_1, \ldots, g_{r+1}) & = & (g_2, \ldots, g_{r+1}) \\  & \vdots &  \\ d_i(g_1, \ldots, g_{r+1}) & = & (g_1, \ldots, g_{i-1}, g_ig_{i+1}, g_{i+2}, \ldots, g_{r+1}) \\  & \vdots &  \\ d_{r+1}(g_1, \ldots, g_{r+1}) & = & (g_1, \ldots, g_{r})
\end{array}
$$
and degeneracy arrows $s_i : G^{r-1} \to G^{r}$ defined by
$$
s_i(g_1, \ldots, g_{r-1}) = (g_1, \ldots, g_i, 1, g_{i+1}, \ldots, g_{r-1}).
$$

\subsection{Lemma} \label{themu} \emph{If $G$ is a group in a category with products, the maps
$$
G^{r+1} \to G^r : (g_1, \ldots, g_{r+1}) \mapsto (g_1^{-1}g_2, g_2^{-1}g_3, \ldots, g_r^{-1}g_{r+1})
$$
define a morphism of simplicial complexes $\nu(\bullet) : G^{prod}(\bullet) \to G^{gr}(\bullet)$.}

\bigskip

{\bf Proof : } This is easily checked.

\bigskip

Let $G$ be a group scheme over a scheme $S$, $p_G : G \to S$ its structural morphism, $1_G : S \hookrightarrow G$ its unit section and
$$
\mu, p_1, p_2 : G \times_S G \to G.
$$
the group law and the projections.

\subsection{Definition} \label{conormal} The \emph{conormal complex of level $m$} of $G$ is the Berthelot-Lieberman filtered complex $\omega^\bullet_{Gm}$ associated to the unit embedding of $S$ into $G^{gr}(\bullet)$.

\bigskip

It will be convenient to write $P_{1m}(r)$ for the $m$-PD-envelope of the unit section in $G^{r}$ and modify the other notations accordingly.

\subsection{Remark} As in the case of the de Rham complex, for which the results are recalled in remark \ref{local}, we have a very simple local description of the situation in the case of a smooth group scheme.
If $s_1, \ldots, s_n$ is a regular sequence of local parameters for the unit section,  then $\mathcal P_{1m}(r)$ is a free module on generators $s^{\{J_1\}} \otimes \cdots \otimes s^{\{J_r\}}$ with $|J_i| \geq 0$ and  $\mathcal I_{1m}^{\{k\}}(r)$ is generated by the $s^{\{J_1\}} \otimes \cdots \otimes s^{\{J_1\}}$ with $|J_1| + \cdots |J_r| \geq k$.
Thus, if $\bar s_i$ denotes the image of $s_i$ in $\omega^1_{Gm}$, we see that $\omega^1_{Gm}$ is a free module on generators $\bar s^J$ with $0 < |J| \leq p^m$ and that, for bigger $r$, $\omega^r_{Gm}$ is generated by the $\bar s^{J_1} \otimes \cdots \otimes \bar s^{J_r}$ with $0 < |J_i| \leq p^m$, subject to some relations.
Of course, the $k$-th step of the Hodge filtration $\Fil^k\omega^r_{Gm}$ has the same generators subject to the additional condition that $|J_1| + \cdots |J_r| \geq k$. In particular, $\Fil^r\omega^r_{Gm} = \omega^r_{Gm}$.

\subsection{Proposition} \emph{There exist canonical morphisms of filtered complexes
$$
\nu(\bullet) : \mathcal P_{1m}(\bullet) \to p_{G*}\mathcal P_{Gm}(\bullet)
$$
and
$$
\nu^\bullet : \omega^\bullet_{Gm} \to p_{G*}\Omega^\bullet_{Gm}.
$$}

{\bf Proof : } Using lemma \ref{themu}, this follows from the functoriality of the construction of the Berthelot-Lieberman complex.

\bigskip

\subsection{Remark} Of course, we also have for each $r$ morphisms of filtered modules
$$
\nu(r) : p_G^*\mathcal P_{1m}(r) \to \mathcal P_{Gm}(r)
$$
and
$$
\nu^r : p_G^*\omega^r_{Gm} \to \Omega^r_{Gm}.
$$

There are also morphisms of filtered complexes
$$
\mathcal P_{1m}(\bullet) = 1_G^*p_G^*\mathcal P_{1m}(\bullet) \to 1_G^*p_G^*p_{G*}\mathcal P_{Gm}(\bullet) \to 1_G^*\mathcal P_{Gm}(\bullet)
$$
and in particular
$$
\omega^\bullet_{Gm} \to 1_G^*\Omega^\bullet_{Gm}.
$$
Without assuming that $G$ is smooth, we cannot really say more.

\bigskip

If $g \in G(S)$, we denote by $T_g : G \to G$ the left translation map given by $h \mapsto gh$.

\subsection{Definition}
A section $\xi$ of $\mathcal P_{Gm}(r)$, or $\mathcal P^k_{Gm}(r)$ is \emph{invariant by translation by $g \in G(S)$} if $T_g^*(\xi) = \xi$.
Moreover, $\xi$ is \emph{invariant by translation} if it is invariant by translation by any section of $G_{S'}$ for any extension $S' \to S$ of the basis.

\bigskip

We indicate with a superscript ${\bullet}^{inv}$ the subsheaf of translation invariant sections.

\subsection{Proposition} \emph{If $G$ is smooth, we have for all $r$ a canonical isomorphism of filtered modules
$$
p_G^*\mathcal P_{1m}(r) \simeq \mathcal P_{Gm}(r)
$$
and therefore also
$$
p_G^*\omega_{Gm}^r \simeq \Omega_{Gm}^r.
$$
We also have canonical  isomorphisms of filtered complexes
$$
\mathcal P_{1m}(\bullet)  \simeq 1_G^*\mathcal P_{Gm}(\bullet) \simeq \mathcal P^{inv}_{Gm}(\bullet)
$$
and therefore also
$$
\omega^\bullet_{Gm}  \simeq 1_G^*\Omega^\bullet_{Gm} \simeq \Omega^{\bullet inv}_{Gm}.
$$}

{\bf Proof : } Since $G$ is smooth, it follows from the second assertion of proposition \ref{gflat} that, for all $r$, we have a canonical isomorphism of filtered modules
$$
p_G^*\mathcal P_{1m}(r) \simeq \mathcal P_{Gm}(r).
$$
The functoriality of the Berthelot-Lieberman construction provides us with a canonical isomorphism of filtered complexes
$$
\mathcal P_{1m}(\bullet) \simeq 1_G^*\mathcal P_{Gm}(\bullet).
$$
We also have injective maps
$$
\mathcal P_{1m}(r) \hookrightarrow p_{G*} \mathcal P_{Gm}(r).
$$
Moreover, since for $g \in G(S)$ we have
$$
\nu(r) \circ (T_g \times \cdots \times T_g)  = \nu(r) : G^{r+1} \longrightarrow G^{r},
$$
we see that the image of $\mathcal P_{1m}(\bullet)$ in $p_{G*}\mathcal P_{Gm}(\bullet)$ is contained in $\mathcal P^{inv}_{Gm}(\bullet)$.
We want to show that the map
$$
\mathcal P_{1m}(\bullet)  \hookrightarrow \mathcal P^{inv}_{Gm}(\bullet)
$$
is bijective.
Since
$$
p_{G*}\mathcal P_{Gm}(r) \simeq p_{G*}\mathcal O_G \otimes_{\mathcal O_S} \mathcal P_{1m}(r)
$$
and $G$ is flat, we are reduced to showing that the canonical map $\mathcal O_S \to \mathcal O_G^{inv}$ is bijective.
If we consider translation by the the identity
$$
G \times G \to G \times G, (g, h) \to (g, gh),$$
we see that, if $f$ is a section of $\mathcal O_G^{inv}$, then $\mu^*(f) = p_1^*(f)$.
Pulling back by $(1_G \times Id_G)^*$, we get $f = p_G^*1_G^*(f)$ which is what we want.

\subsection{Proposition}  \emph{If $G$ is smooth and if we let
$$\delta := p^*_2 - \mu^* + p^*_1 : p_{G*}\mathcal P_{Gm}(1)  \to p_{G^2*}\mathcal P_{G^2m}(1),
$$
we have an isomorphism of filtered modules
$$
p_{G*}\mathcal I_{Gm}(1) \cap \ker \delta = \mathcal I^{inv}_{Gm}(1) \cap \ker d
$$
and therefore also
$$
p_{G*}\Omega^1_{Gm} \cap \ker \delta = \Omega^{1inv}_{Gm} \cap \ker d.
$$}
\bigskip

{\bf Proof : } If $g \in G(S)$ and $i_g  : G \to G \times G, h \mapsto (g, h)$, we have $p_2 \circ i_g = \Id_G, \mu \circ i_g = T_g$ and $p_1 \circ i_g = g \circ p_G$.
It follows that
$$
i_g^* \circ \delta = (\Id - T^*_g) + p_G^* \circ g^*.
$$
If $\xi \in p_{G*}\mathcal I_{Gm}(1)$ then, trivially,  $g^*(\xi) = 0$ and it follows that $i_g^*(\delta(\xi)) = \xi - T^*_g(\xi)$.
Thus, any $\xi \in p_{G*}\mathcal I_{Gm}(1) \cap \ker \delta$ is $g$-invariant.
This is true for any $g$ and after any base change.
It follows that
$$
p_{G*}\mathcal I_{Gm}(1) \cap \ker \delta \subset \mathcal I^{inv}_{Gm}(1).
$$
By functoriality, there is a commutative diagram

$$
\begin{array}{ccc }
\mathcal P_{1m}(1) & \stackrel d\longrightarrow & \mathcal P_{1m}(2) \\
\downarrow & & \downarrow \\
p_{G*}\mathcal P_{Gm}(1) & \stackrel \delta \longrightarrow & p_{G^2*}\mathcal P_{G^2m}(1)
\end{array}
$$

Since we assumed that $G$ is smooth, the vertical arrows induce filtered isomorphisms with the invariant part in the bottom and we obtain

$$
\begin{array}{ccc }
\mathcal P_{1m}(1) & \stackrel d\longrightarrow & \mathcal P_{1m}(2) \\
|| & & || \\
\mathcal P^{inv}_{Gm}(1) & \stackrel \delta \longrightarrow & \mathcal P^{inv}_{G^2m}(1)
\end{array}
$$

from which it follows that $d = \delta$ on $P^{inv}_{Gm}(1)$ and the result follows.

\subsection{Notations} \label {invariant} The sheaf
$$
p_{G*}\Omega^1_{Gm} \cap \ker \delta = \Omega^{1inv}_{Gm} \cap \ker d \simeq \omega^1_{Gm} \cap \ker d \simeq \mathcal H^1(\omega^\bullet_{Gm})
$$
of closed invariant forms will be generally written $\omega^{(m)}_{G,0}$.

\bigskip

Note that in the case $m = 0$, the situation was a lot simpler since invariant forms are automatically closed, and in particular, the kernel of $\delta$ is exactly the sheaf of invariant differential forms.

\subsection{Examples}

\begin{enumerate}

\item If $G = \mathbf G_{aS}$ with parameter $t$, we already know that $\omega^1_{m}$ is the free module on $\bar t, \bar t^2, \ldots, \bar t^{p^m}$ and it is not difficult to verify that $\bar t \in \omega^{(m)}_{0}$.
More precisely, the group law is given by $t \mapsto t \otimes 1 + 1 \otimes t$ and it follows that
$$
\delta(t^k) = t^k \otimes 1 - (t \otimes 1 + 1 \otimes t)^k + 1 \otimes t^k = \sum_{i=1}^{k-1} {k \choose i} t^i \otimes t^{k-i}.
$$
But we also see that, if $\Char S = p$, then $\bar t^{p^j} \in \omega^{(m)}_{0}$ for all $j$.
Actually, in this case, $\omega^{(m)}_{0}$ is locally free on the generators $\bar t, \bar t^p, \ldots \bar t^{p^m}$.
In particular, we see that the filtration on $\omega^{(m)}_{0}$ can have length exactly $p^m - 1$.

In the case $S = \Spec \Z/4$ (so that $p = 2$) and $m = 1$, we see that $\delta(\bar t^2) = 2 \bar t \otimes \bar t \not = 0$ but that $2\bar t^2 \in \omega^{(1)}_{0}$.
Thus, $\omega^{(m)}_{0}$ is not always locally free.

The opposite map in $\mathbf G_a$ is given by $t \mapsto -t$ and the difference is therefore given by $t \mapsto - t \otimes 1 + 1 \otimes t$.
It follows that the canonical inclusion $\omega^1_{m} \hookrightarrow p_*\Omega^1_{m}$ sends $\bar t$ to $dt$ and therefore, $dt$ is an invariant differential of level $m$ on $\mathbf G_{aS}$.

\item If $G = \mathbf G_{mS}$ with parameter $t$, then $\omega^1_{m}$ is the free module on $\bar s, \bar s^2, \ldots, \bar s^{p^m}$ with $s = t -1$.
Let us consider $\log t^{p^m} := \log(1 + s)^{p^m} \in \widehat{\mathcal P_{1m}(1)}$.
Using the fact that the group law is given by $t \mapsto t \otimes t$, it is not difficult to check that $\delta(\log t^{p^m}) = 0$.
Actually, this is a purely formal calculation that can be done over $\C$ where this is well known.
It follows that
$$
\sum_{i=1}^{p^m} \frac {p^{m}}i \bar s^i \in \omega_0^{(m)}.
$$
Note that, unlike the case $m = 0$, $\Fil^2 \omega^{(m)}_{0} \not = 0$ in general.
For example, if $p = 2$, $m = 1$ and $S = \Spec \Z/4$, then we have $2\bar s - \bar s^2 \in \Fil^1$ and $2\bar s^2 \in \Fil^2$.

Since the inverse on $\mathbf G_m$ is given by $t \mapsto t^{-1}$, the difference is given by $t \mapsto t^{-1} \otimes t$.
Thus we see that $t - 1$ is sent to $t^{-1} \otimes t  - 1 \otimes 1 = t^{-1} (1 \otimes t  - t \otimes 1)$ and it follows that the canonical inclusion $\omega^1_{m} \hookrightarrow p_*\Omega^1_{m}$ sends $\bar s$ to $dt/t$.
Finally, we get that
$$
\sum_{i=1}^{p^m} \frac {p^{m}}i (\frac {dt} t)^i
$$
is an invariant differential of level $m$ on $\mathbf G_{mS}$.

\item We assume now $p \not = 2$ and consider the Legendre elliptic curve $E$ given by the equation
$$
y^2 = x(x - 1)(x - \lambda)
$$
over $\mathbf A^1_S \backslash \{0, 1\}$.
Since we are interested in the behavior at $O$ which is the point at infinity, we make the usual change of coordinates $z = -x/y, w = -1/y$ and the equation becomes
$$
w = z(z - w)(z - \lambda w).
$$
Thus, $P_{1m}(1)$ is the $m$-th divided power envelope of $\mathcal O_{\mathbf A^1_S \backslash \{0, 1\}}[z]$ with respect to $z$ and $\omega^1_{Em}$ is the free module on $\bar z, \ldots, \bar z^{p^m}$.
A quick calculation as explained in Chapter 4 of \cite {Silverman86}, for example, shows that the group law is given by
$$
z \mapsto 1 \otimes z + z \otimes 1 + (\lambda + 1)(z \otimes z^2 + z^2 \otimes z) + \cdots
$$
Assume now that $p = 3$ and $m = 1$.
Taking into account the symmetry of the above series due to the commutativity of the group law, one sees that higher powers don't play any role and that
$$
\delta (\bar z) = (\lambda + 1)(\bar z \otimes \bar z^2 + \bar z^2 \otimes \bar z),
$$
$$
\delta(\bar z^2) = 2\bar z \otimes \bar z + (\lambda + 1)(\bar z^2 \otimes \bar z^2),
$$
$$
\delta(\bar z^3) = -3(\bar z \otimes \bar z^2 + \bar z^2 \otimes \bar z).
$$
It turns out that $3 \bar  z + (1 + \lambda) \bar z^3$ is in $\omega^{(1)}_{E,0}$.

\bigskip

Assume moreover that $\Char S = 3$.
Then we get $(1 + \lambda)\bar z^3$ which is in $\Fil^3$.
It follows from the formulas that $\omega^{(1)}_{E,0}$ is actually free of rank 1 with generator $\bar z^3$ outside \linebreak $\lambda = -1$.
However, at the fiber $\lambda = -1$, both $\bar z$ and $\bar z^3$ are in $\omega^{(1)}_{E,0}$.
In other words, the supersingular fiber is characterized by the fact that $\omega^{(1)}_{E,0}$ is of rank $2$ in contrast with the general case where it has rank $1$.
This is a phenomenon that is specific to higher level because in the classical situation, $\omega^{(0)}_{E,0} = \omega_E$ is globally free of rank 1. Note however that in the case $m = 1$ and $p = 3$, then $Fil^{p^m}\omega^{(m)}_{E,0}$ is also globally free of rank $1$.

The next question would be to describe the canonical inclusion $\omega^1_m \hookrightarrow p_*\Omega^1_m$.
Since taking opposite on $E$ is given by horizontal symmetry, it sends $z$ to $-z$, and it follows that the difference in $E$ is given by
$$
z \mapsto 1 \otimes z - z \otimes 1 + (\lambda + 1)(- z \otimes z^2 + z^2 \otimes z) + \cdots
$$
When $m = 0$, one can check that $\bar z$ is sent to $dx/2y$ but we have been to lazy to do the computations in other cases.

\end{enumerate}

\subsection{Proposition} \emph{If $A$ is an abelian scheme over $S$, there is a Hodge spectral sequence
$$
E_2^{ij} = R^jp_{A*}\mathcal O_A \otimes_{\mathcal O_S} \mathcal H^i(\omega_{Am}^\bullet) \Rightarrow \mathcal H^n_{dRm}(A/S)
$$}

{\bf Proof : } Consider the Hodge to de Rham spectral sequence of level $m$ :
$$
E_1^{ij} = R^jp_{A*}\Omega_{Am}^i \Rightarrow \mathcal H^n_{dRm}(A/S)
$$
We know that, for all $i$, $\Omega_{Am}^i  = p_A^*\omega_{Am}^i$.
Since we also have $p_{A*}\mathcal O_A = \mathcal O_S$, we see that
$$
R^jp_{A*}\Omega_{Am}^\bullet  = R^jp_{A*}\mathcal O_A \otimes_{\mathcal O_S} \omega_{Am}^\bullet
$$
and we get the spectral sequence
$$
E_2^{ij} = R^jp_{A*}\mathcal O_A \otimes_{\mathcal O_S} \mathcal Hä(\omega_{Am}^\bullet) \Rightarrow \mathcal H^n_{dRm}(A/S).
$$

\subsection{Corollary} \emph{We have an exact sequence
$$
0 \to \mathcal \omega^{(m)}_{A,0} \to \mathcal H^1_{dRm}(A/S) \to \mathcal R^1p_{A*}\mathcal O_A
$$}

\bigskip

{\bf Proof : } Just note that $\omega^{(m)}_{A,0} = H^1(\omega^\bullet_{Am})$.

\subsection{Remark} When $S$ is smooth, using some results from \cite {Trihan98}, we can show that we actually get a short exact sequence.
We hope to be able to prove this in general in a forthcoming article.

%
%

\section{Crystalline extension groups}

In this section, we generalize to higher level the first results of the second chapter of \cite{BerthelotBreenMessing82}.

\subsection{Definitions} A category $I$ is said \emph{very small} if $\Ob(I)$ is countable and for each $\alpha, \beta \in \Ob(I)$, $\Hom(\alpha, \beta)$ is finite.
A \emph{decoration} on a very small category $I$ is a pair composed by a \emph{degree map} $d : \Ob(I) \to \Z$ with finite fibers and a \emph{sign map} $\epsilon : \Arr(I) \to \{\pm\}$.

\subsection{Notations} If $Y$ is any object in a topos $\mathcal T$, we will denote by $\Z^{(Y)}$ the free abelian group on $Y$ in $\mathcal T$.

\bigskip

Let $(I, d, \epsilon)$ be a decorated category and $Y_\bullet := (Y_\alpha, f_\lambda)$ a diagram indexed by $I$ in $\mathcal T$.
We set
$$
C_n(Y_\bullet) := \oplus_{d(\alpha) = n} \Z^{(Y_{\alpha})}
$$
and let $\delta_n : C_n \to C_{n-1}$ be given by
$$
\sum_{\lambda : \alpha \to \beta} \epsilon (\lambda) f_\lambda : \Z^{(Y_{\alpha})} \to \Z^{(Y_{\beta})}
$$
whenever $\deg(\alpha) = n$ and $\deg(\beta) = n - 1$.

\subsection{Definition} A diagram $Y_\bullet$ indexed by a decorated category $(I, d, \epsilon)$ is \emph{nice} if $C(Y_\bullet)$ is a complex.

\subsection{Remark} \label {groupalso} P. Deligne has shown that any abelian group $G$ in a topos $\mathcal T$ has a canonical left resolution $C_{\bullet}(G)$ that fits in the above setting.
More precisely, there exists a decorated category $I$ and a nice diagram $\Delta(G)_{\bullet} := (G^{n_\alpha}, f_\lambda)$ indexed by $I$ in the subcategory of $\mathcal T$ generated by the powers of $G$, the projections and the group law such that $C_{\bullet}(G) = C_{\bullet}(\Delta(G)_{\bullet})$.

\bigskip

To the best of our knowledge, this result of Deligne is unpublished, and we can offer no better reference than what is said in \cite {BerthelotBreenMessing82}.

\bigskip

On the other hand, we are mainly interested in the first terms of the complex $C_{\bullet}(G)$ in which case we may use the description given in section 2.1 of \cite {BerthelotBreenMessing82} which we now recall for the reader's convenience.
The diagram has only one object $G$ in degree $0$.
In degree 1, there is also a unique object $G^2$ and the morphisms $G^2 \to G$ are

$$
\begin{array}{cccc}G^2 & \to & G & \textrm{sign} \\ (g_1, g_2) & \mapsto & g_1 & - \\(g_1, g_2) & \mapsto & g_1 + g_2 & + \\ (g_1, g_2) & \mapsto & g_2 & -\end{array}
$$

In degree $3$, there are two objects, namely $G^2$ and $G^3$.
The morphisms from degree 3 to degree 2 are given by
$$
\begin{array}{cccc}G^2 & \to & G^2 & \textrm{sign} \\(g_1, g_2) & \to & (g_1, g_2) & + \\(g_1, g_2) & \to & (g_2, g_1) & -\end{array}
$$
and by
$$
\begin{array}{cccc}G^3 & \to & G^2 & \textrm{sign} \\(g_1, g_2, g_3) & \mapsto & (g_2, g_3) & - \\(g_1, g_2, g_3) & \mapsto & (g_1 + g_2, g_3) & + \\(g_1, g_2, g_3) & \mapsto & (g_1, g_2 + g_3) & - \\(g_1, g_2, g^3) & \mapsto & (g_1, g_2) & +\end{array}
$$
Thus  we get
$$
C_{\bullet}(G) = \cdots \to \Z^{(G^3)} \oplus \Z^{(G^2)} \to \Z^{(G^2)} \to \Z^{(G)}
$$
and the maps are given by
$$
\begin{array}{ccl}[g_1, g_2] & \mapsto & - [g_1] + [g_1 + g_2] - [g_2]\end{array}
$$
from degree 1 to degree 0 and
$$
\begin{array}{ccl}
[g_1, g_2] & \mapsto & [g_1, g_2] - [g_2, g_1]
\\ \lbrack g_1, g_2, g_3] & \mapsto & - [g_2, g_3] + [g_1 + g_2, g_3] - [g_1, g_2 + g_3] + [g_1, g_2]
.\end{array}
$$
from degree 2 to degree 1.

\bigskip

\subsection{Lemma} \label {specseq} \emph{Let $Y_\bullet$ be a nice diagram indexed by a decorated category $(I, d, \epsilon)$ in some topos $\mathcal T$ and $E$ an abelian sheaf in $\mathcal T$.
Then, there is a spectral sequence
$$
E_1^{r,s} = \bigoplus_{d(\alpha) = r} R^sj_{Y_\alpha*}j_{Y_\alpha}^{-1}E \Rightarrow \mathcal Ext^{r+s}(C(Y_\bullet), E).
$$}

\bigskip

{\bf Proof : } If $Y$ is any object of $\mathcal T$, we have
$$
\mathcal Hom_{Gr}(\Z^{(Y)}, E) = \mathcal Hom(Y, E) = j_{Y*}j_Y^{-1} E.
$$
It follows that the complex $\mathcal Hom_{Gr}(C_\bullet(Y_\bullet), E)$ is canonically isomorphic to a complex whose terms are all of the form $\oplus_{\deg \alpha = n}  j_{Y_\alpha*}j_{Y_\alpha}^{-1}E$.
Since pulling back by localization is exact and preserves injective sheaves, we get the spectral sequence by applying this remark to an injective resolution of $E$.

\bigskip

Let $(S, \mathfrak a, \mathfrak b)$ be an $m$-PD-scheme with $p$ locally nilpotent and $p \in \mathfrak a$.

\bigskip

Let $X$ be an $S$-scheme to which the $m$-PD-structure of $S$ extends.

\subsection{Remark} We denote by $\Sch'_{/X}$ the category of $X$-schemes to which the $m$-PD-structure of $S$ extends.
The Zariski topology on this category is coarser than the canonical topology and we obtain an embedding of $\Sch'_{/X}$ into the corresponding topos $X_{\ZAR'}$.
We will now use the big crystalline topos of level $m$, $(X/S)^{(m)}_{\CRIS}$ that was introduced in \cite{EtesseLeStum97} (but considering only schemes to which the $m$-PD structure extends).
Composing the embedding $\Sch'_{/X} \hookrightarrow X_{\ZAR'}$ with the canonical map
$$
v_{X/S*} : X_{\ZAR'} \to  (X/S)^{(m)}_{\CRIS}
$$
from section 1.10 of \cite{EtesseLeStum97} gives a functor
$$
\begin{array}{ccc}\Sch'_{/X} & \to & (X/S)^{(m)}_{\CRIS} \\Y & \mapsto & \underline Y\end{array}
$$
Note that if $Y \in \Sch'_{/X}$ with structural morphism $f_Y : Y \to X$ then the canonical morphism
$$
f_{Y/X} : (Y/S)^{(m)}_{\CRIS} \to (X/S)^{(m)}_{\CRIS}
$$
factors as an isomorphism $(Y/S)^{(m)}_{\CRIS} \simeq (X/S)^{(m)}_{\CRIS/\underline Y}$ followed by the localization map
$$
j_{\underline Y} : (X/S)^{(m)}_{\CRIS/\underline Y} \to (X/S)^{(m)}_{\CRIS}.
$$

\subsection{Proposition} \label{bigcris} \emph{Let $Y_\bullet$ be a nice diagram in $\Sch'_{/X}$ indexed by a decorated category $(I, d, \epsilon)$.
If $E$ is an abelian sheaf on $\CRIS^{(m)}(X/S)$, there is a spectral sequence
$$
E_1^{r,s} = \bigoplus_{d(\alpha) = r} R^sf_{Y_\alpha/X*}f_{Y_\alpha/X}^{-1} E \Rightarrow \mathcal Ext^{r+s}(C(\underline Y_\bullet), E).
$$}

\bigskip

{\bf Proof : } Taking into account the previous remark, this immediately follows from lemma \ref{specseq}.

\subsection{Remarks} Just as in Chapter III, section 4 of \cite {Berthelot74}, if $Y \in \Sch'_{X}$, there is a morphism of topos
$$
(Y/S)^{(m)}_{\cris} \to (X/S)^{(m)}_{\CRIS}
$$
whose inverse image functor $E \mapsto E_Y$ might be called restriction.
For any sheaf $E$ on $\CRIS^{(m)}(X/S)$ and any morphism $g :Y' \to Y$, there is a canonical transition map \linebreak $g^{-1}E_Y \to E_{Y'}$ and these data uniquely determine $E$.
An $m$-crystal on the big site can be defined as an $\mathcal O^{(m)}_{X/S}$-module $E$ such that all $E_Y$ are $m$-crystals and the transition maps induce isomorphisms $g^{*}E_Y \simeq E_{Y'}$.
In particular, the functor $E \mapsto E_X$ is an equivalence of categories between $m$-crystals on $\CRIS^{(m)}(X/S)$ and $m$-crystals on $\Cris^{(m)}(X/S)$.
Finally, note that any filtered (resp. transversal) $\mathcal O^{(m)}_{X/S}$-module $(E, \Fil^\bullet)$ on $\CRIS^{(m)}(X/S)$ restricts for each $Y$ to a filtered (resp. transversal) $\mathcal O^{(m)}_{Y/S}$-module $(E_Y, \Fil^\bullet E_Y)$ on $\Cris^{(m)}(Y/S)$.

\bigskip

For future reference, note also that, as for level 0, in which case this is proved in 1.1.16.4 of \cite {BerthelotBreenMessing82}, we have for any abelian sheaf on $\CRIS^{(m)}(Y/S)$,
$$
(Rf_{Y/X*} E)_{(U,T)} = Rf_{U \times Y/T*}E_{U \times Y}.
$$

\subsection{Definition} A \emph{big transversal $m$-crystal} on $X/S$ is a crystal $E$ on $\CRIS^{(m)}(X/S)$, endowed with a filtration $\Fil^\bullet$ such that for each $Y \in \Sch'_{/X}$, the filtered $m$-crystal $(E_Y, \Fil^\bullet E_Y)$ is a transversal $m$-crystal.

\subsection{Lemma} \emph{Let $Y \in \Sch'_{/X}$ and $(U \hookrightarrow T) \in \CRIS^{(m)}(X/S)$.
Let $i : U \times_X Y \hookrightarrow Z$ be an immersion into a smooth $Z$ in $\Sch'_{/T}$.
If $E$ is a big transversal $m$-crystal on $X/S$, we have a canonical isomorphism
$$
(R^sf_{Y/X*}f_{Y/X}^{-1} \Fil^kE)_{(U,T)} = R^sf_{Z*}\Fil^k[(i_{\cris*}E_{U \times Y})_{Z} \otimes \Omega^\bullet_{Z/Tm}].
$$}

\bigskip

{\bf Proof : } We already know that
$$
(R^sf_{Y/X*}f_{Y/X}^{-1} \Fil^k E)_{(U,T)} = R^sf_{U \times Y/T*}\Fil^kE_{U \times Y}.
$$
On the other hand, the filtered Poincar\'e lemma in level $m$ tells us that
$$
\Fil^ki_*Ru_{Y*} E_{U \times Y} \simeq \Fil^k[(i_{\cris*}E_{U \times Y})_Z \otimes \Omega^\bullet_{Z/Tm}]
$$
Applying $R^sf_{Z*}$ gives
$$R^sf_{U \times Y/T*} \Fil^k E_{U \times Y} \simeq Rf^s_{Z*}\Fil^k[(i_{\cris*}E_{U \times Y})_Z \otimes \Omega^\bullet_{Z/Tm}]
$$
and we are done.

\subsection{Proposition} \label{SpecSeq}  \emph{Let $Y_\bullet$ be a nice diagram in $\Sch'_{/X}$ and $(U \hookrightarrow T) \in \CRIS^{(m)}(X/S)$. Let $Z_\bullet$ be a nice diagram in $\Sch'_{/T}$ with all $Z_\alpha$ smooth and $i_\bullet : U \times_X Y_\bullet \hookrightarrow Z_\bullet$ a compatible family of immersions.
If $E$ is a big transversal $m$-crystal on $X/S$, we have a canonical spectral sequence
$$
E_1^{r,s} = \bigoplus_{d(\alpha) = r} R^sf_{Z_\alpha*}\Fil^k[(i_{\alpha \cris*}E_{U \times Y_\alpha})_{Z_{\alpha}} \otimes \Omega^\bullet_{Z_{\alpha}/Tm}] \Rightarrow \mathcal Ext^{r+s}(C(\underline Y_\bullet), \Fil^k E)_{(U,T)}.
$$}

\bigskip

{\bf Proof : } Since taking value on some object $(U, T)$ is an exact functor, we know from proposition \ref{bigcris} that there is a spectral sequence
$$
E_1^{r,s} = \bigoplus_{d(\alpha) = r} [R^sf_{Y_\alpha/X*}f_{Y_\alpha/X}^{-1} \Fil^kE]_{(U,T)}
\Rightarrow \mathcal Ext^{r+s}(C(\underline Y_\bullet), \Fil^kE)_{(U,T)}.
$$
It is therefore sufficient to prove that for each $\alpha$, we have
$$
[R^sf_{Y_\alpha/X*}f_{Y_\alpha/X}^{-1} \Fil^kE]_{(U,T)} = R^sf_{Z_\alpha*}\Fil^k[(i_{\alpha \cris*}E_{U \times Y_\alpha})_{Z_{\alpha}} \otimes \Omega^\bullet_{Z_{\alpha}/Tm}].
$$
and this follows from the lemma.

\bigskip

Applying this to the nice diagram $\Delta(G)_{\bullet}$ that, as we remarked in \ref{groupalso}, gives rise to Deligne's resolution, we immediately get

\subsection{Theorem} \label{SpecGroup} \emph{Let $G$ be a group scheme in $\Sch'_{/X}$ and $(U \hookrightarrow T) \in \CRIS^{(m)}(X/S)$.
Let $H$ be a smooth group scheme on $T$ and $i : G_U \hookrightarrow H_U$ an immersion of groups.
If $E$ is a big transversal $m$-crystal on $X/S$, we have a canonical spectral sequence
$$
E_1^{r,s} = \bigoplus_{d(\alpha) = r} R^sf_{H^{n_\alpha}*}\Fil^k[(i_{\alpha \cris*}E_{G_U^{n_\alpha}})_{H^{n_\alpha}} \otimes \Omega^\bullet_{H^{n_\alpha}/Tm}] \Rightarrow \mathcal Ext^{r+s}(\underline G, \Fil^k E)_{(U,T)}.
$$}

\subsection{Remarks} \label{complex} If $E$ is any $m$-crystal on $X/S$, then the trivial effective filtration $\mathcal I^{\{k\}}_{X/S}E$ turns it into a big transversal $m$-crystal and both proposition \ref{SpecSeq} and theorem \ref{SpecGroup} apply.
Also, the theorem is still valid if we replace $G$ with a complex of abelian groups.

\bigskip

Note also that, when $H$ is affine, we have, as in theorem 2.1.8 of \cite{BerthelotBreenMessing82}, an isomorphism in the derived category
$$
\Fil^k [\bigoplus_{d(\alpha) = \bullet} f_{H^{n_\alpha}*}[(i_{\alpha \cris*}E_{G_U^{n_\alpha}})_{H^{n_\alpha}} \otimes \Omega^\bullet_{H^{n_\alpha}/Tm}]_s \simeq R\mathcal Hom(\underline G, \Fil^k E)_{(U,T)}.
$$
The proof goes exactly as in \cite{BerthelotBreenMessing82}.

\subsection{Proposition} \label{exti} \emph{Let $G$ be a smooth group scheme in $\Sch'_{/X}$ and $U  \in \Sch_{/X}$.
If $k > 0$, we have
$$
\begin{array}{ccc} \mathcal Hom(\underline G, \mathcal I^{\{k\}}_{X/S})_{(U,U)} & = & 0 \\\mathcal Ext^1(\underline G,  \mathcal I^{\{k\}}_{X/S})_{(U,U)} & = & Fil^k \omega^{(m)}_{G_U,0}\end{array}
$$
where $\omega^{(m)}_{G_U,0}$ denotes, as in \ref{invariant}, the sheaf of closed invariant forms of level $m$ on $G_U$.
In particular,
$$\mathcal Ext^1(\underline G,  \mathcal I_{X/S})_{(U,U)} = \omega^{(m)}_{G_U,0}.
$$
}

\bigskip

{\bf Proof : } We consider the spectral sequence
$$
E_1^{r,s} = \bigoplus_{d(\alpha) = r} R^sf_{G_U^{n_\alpha}*}\Fil^k \Omega^\bullet_{G_U^{n_\alpha}/Um} \Rightarrow \mathcal Ext^{p+q}(\underline G, \mathcal I^{\{k\}}_{X/S})_{(U,U)}.
$$
Note first that $E_1^{r,s} = 0$ for $r < 0$ or $s < 0$.
Actually, since $k > 0$, we have $\Fil^k \mathcal O_{G_U^{n_\alpha}} = 0$ and it follows that $E_1^{r,s} = 0$ for $s = 0$ also.
Therefore
$$
\mathcal Hom(\underline G, \mathcal I^{\{k\}}_{X/S})_{(U,U)}  =  0
$$
and
$$
\mathcal Ext^1(\underline G,  \mathcal I^{\{k\}}_{X/S})_{(U,U)} = E_2^{0,1} = \ker \delta : E_1^{0,1} \to E_1^{1,1}
$$
But we have for all $\alpha$,
$$
R^1f_{G_U^{n_\alpha}*}\Fil^k \Omega^\bullet_{G_U^{n_\alpha}/Um} = \Fil^k \ker [d : f_{G_U^{n_\alpha}*}\Omega^1_{G_U^{n_\alpha}/Um}  \to f_{G_U^{n_\alpha}*}\Omega^2_{G_U^{n_\alpha}/Um}]
$$
and it follows that
$$E_2^{0,1} = \Fil^k \ker [d : f_{G_U*}\Omega^1_{G_U/Um} \to f_{G_U*}\Omega^2_{G_U^/Um}] \cap \ker[\delta : f_{G_U*}\Omega^1_{G_U/Um} \to f_{G_U^2*}\Omega^1_{G_U^2/Um}]
$$
which is exactly $Fil^k \omega^{(m)}_{G_U,0}$.

\subsection{Remark} The proofs we have presented for theorem \ref{SpecGroup} and proposition \ref{exti} rely on the unpublished theorem of Deligne mentioned in \ref{groupalso}.
We have taken this option because we think it is the most elegant and natural (even more so if one wants to prove the more general form of theorem \ref{SpecGroup} suggested in \ref{complex}).
For readers who fill uncomfortable using unpublished results, we should mention that it would not be difficult to avoid Deligne's resolution and work with a partial resolution as in \cite{BerthelotBreenMessing82}, chapter 2, section 1.

%
%


\bibliographystyle{plain}

%
%

\vskip 12pt {\sc IRMAR, Universit\'e de Rennes I, Campus de Beaulieu, F-35042 Rennes Cedex, France\\ e-mail:}
bernard.le-stum@univ-rennes1.fr

\vskip 12pt {\sc Departamento de Matem\'aticas, Universidad Aut\'onoma de Madrid, E-28049 Madrid, Spain \\ e-mail:}
adolfo.quiros@uam.es

\end{document}